\newcommand{\R}{\mathbb{R}}

\newcommand{\A}{\mathscr{A}}
\newcommand{\E}{\mathscr{E}}
\newcommand{\m}{\mathfrak{m}}

\newcommand{\Diff}{{\rm Diff}}   
\mathchardef\varepsilon="010F
\mathchardef\epsilon="0122
\mathchardef\vartheta="0112
\mathchardef\theta="0123
\mathchardef\varrho="011A
\mathchardef\rho="0125
\mathchardef\varphi="011E     
\mathchardef\phi="0127
\renewcommand \emptyset \varnothing
\documentclass[12pt]{article}            
\usepackage[all]{xy}
\usepackage[latin1]{inputenc}        
\usepackage[dvips]{graphics}
\usepackage[dvips]{graphicx}
\usepackage{amsfonts}
\usepackage{amssymb}
\usepackage{amsmath}
\usepackage{amstext}
\usepackage{amsbsy}
\usepackage{amsopn}
\usepackage{amsthm}
\usepackage{amscd}
\usepackage{amsxtra}
\usepackage{mathrsfs}
\usepackage{upref}
\author{{\sc Gianmarco Capitanio}\\
{\it \small Universit{\' e} D. Diderot -- Paris VII, 
Equipe de G{\' e}om{\' e}trie et Dynamique } \\
{\it \small Case 7012 -- 2, place Jussieu , 75251 Paris Cedex 05} \\ 
{\small {\it e-mail:} {\rm Gianmarco.Capitanio@math.jussieu.fr}}
}
\date{}                     
\title{Legendrian graphs generated by tangential families}
\begin{document}        
\maketitle
\theoremstyle{plain}
\newtheorem{theorem*}{\bf Theorem}
\newtheorem{theorem}{\bf Theorem}
\newtheorem*{conjecture*}{\bf Conjecture}
\newtheorem{lemma}{\bf Lemma}
\newtheorem{proposition}{\bf Proposition}
\newtheorem*{proposition*}{\bf Proposition}
\newtheorem{corollary}{\bf Corollary}
\newtheorem*{corollary*}{\bf Corollary}
\theoremstyle{definition}
\newtheorem*{definition*}{\bf Definition}
\newtheorem*{definitions*}{\bf Definitions}

\newtheorem{example}{\bf Example}
\newtheorem*{example*}{\bf Example}
\theoremstyle{remark}
\newtheorem*{remark*}{\bf Remark}
\newtheorem{remark}{\bf Remark}
\newtheorem*{remarks*}{\bf Remarks}
\newtheorem{remarks}{\bf Remarks}

\begin{abstract}
We construct a Legendrian version of Envelope theory. 
A tangential family is a $1$-parameter family of rays emanating 
tangentially from a smooth plane curve.
The Legendrian graph of the family is the union of the
Legendrian lifts of the family curves in the projectivized cotangent
bundle $PT^*\R^2$.  
We study the singularities of Legendrian graphs and 
their stability under small tangential deformations.  
We also find normal forms of their projections into the plane. 
This allows to interprete the beaks perestroika as the apparent
contour of a deformation of the Double Whitney Umbrella singularity
$A_1^\pm$.   
\end{abstract}

{\small {\bf \sc Keywords :} Envelope theory,
Tangential families, Contact Geometry.}

{\small {\bf \sc 2000 MSC :} 14B05, 53C15, 58K25.}

\section{Introduction}\label{sct:1}

A tangential family is a $1$-parameter family of ``rays'' emanating tangentially 
from a smooth plane curve. 
Tangential families and their envelopes (or caustics) are natural objects in 
Differential Geometry: for instance, every curve in a Riemannian surface 
defines the tangential family of its tangent geodesics. 
The theory of tangential families is related to the study 
developed by Thom and Arnold for plane envelopes 
(see \cite{arnold1976a}, \cite{arnold1976b} and \cite{Thom}). 
In \cite{stable} and \cite{simple} we studied stable and simple
singularities of tangential family germs 
(with respect to deformations among tangential families). 

In this paper we construct a Legendrian version of tangential family
theory.  
The envelope of a tangential family is viewed as the 
apparent contour of the surface, called Legendrian graph, formed by the
union of the Legendrian lifts of the family curves in the
projectivized cotangent bundle of the plane. 

We classify the Legendrian graph singularities that are stable 
under small tangential deformations of the generating tangential families.  
We prove that, in addition to a regular Legendrian graph, there exists
just one more local stable singularity, the Double
Whitney Umbrella $A_1^\pm$. 
Furthermore, we find normal forms of typical projections of 
Legendrian graphs into the plane. 
This allows to interprete the beaks perestroika as the apparent
contour of a non-tangential deformation of the Double Whitney Umbrella 
singularity in the projectivization of the cotangent bundle $T^*\R^2$.  

Our results are related to several theories, concerning 
Maps from the Space to the Plane (Mond \cite{Mond}), 
Projections of Manifolds with Boundaries (Goryunov
\cite{goryunov2}, Bruce and Giblin \cite{BG}), 
Singular Lagrangian Varieties and their Lagrangian mappings
(Givental \cite{givental1}).  

\noindent {\bf Acknowledgments.} 
This paper contains a part of the results of my PhD Thesis
(\cite{these}). 
I wish to express my deep gratitude to my advisor V.I. Arnold. 

\section{Legendrian graphs and their singularities}\label{sct:2}

Unless otherwise specified, all the objects considered below 
are supposed to be of class $\mathscr C^\infty$; 
by plane curve we mean an embedded $1$-submanifold of the plane.   

In this section we recall basic facts about tangential families 
and we define their Legendrian graphs in the projectivized 
cotangent bundle $PT^*\R^2$. 
We study the typical singularities of these graphs up to Left-Right
equivalence.   
This classification considers neither the fiber  
nor the contact structure of $PT^*\R^2$. 
A classification of Legendrian graphs taking into account 
the fiber bundle structure is the object of Section \ref{sct:3}. 
 
Let $f:\R^2\rightarrow \R^2$ be a mapping of the source plane, 
equipped with the coordinates $\xi$ and $t$, to another plane. 
If $\partial_t f$ vanishes nowhere, then $f_\xi:=f(\xi,\cdot)$
parameterizes an immersed curve $\Gamma_\xi$. 
Hence, $f$ parameterizes the $1$-parameter family of curves
$\{\Gamma_\xi:\xi\in\R\}$. 

\begin{definition*}
The family parameterized by $f$ is a {\it tangential family} 
if $f(\cdot,0)$ parameterizes an embedded smooth curve, 
called the {\it support}, and $\Gamma_\xi$ is tangent to $\gamma$ at
$f(\xi,0)$ for every $\xi\in\R$. 
\end{definition*}

The {\it graph} of the family is the 
surface $\Phi:= \cup_{\xi\in\R} \{(\xi,\Gamma_\xi)\}\subset\R^3$. 
The {\it envelope} is the apparent contour of $\Phi$ under 
the projection $\pi: \R^3\rightarrow \R^2$, $\pi (\xi,P):=P$ (i.e. the
critical value set of $\pi|_\Phi$);   
the {\it criminant set} is the critical set of $\pi|_\Phi$.
By the very definition, the support of a tangential 
family belongs to its envelope. 

A $p$-parameter deformation $F:\R^2\times \R^p\rightarrow \R^2$ 
of a tangential family $f$ is {\it tangential} if 
$F_\lambda:=F(\cdot;\lambda)$ is a tangential family for every 
$\lambda$.  
Remark that the supports of the deformed families form a
smooth deformation of the support of the initial family. 

Below we will consider tangential family germs. 
Note that graphs of tangential family germs are smooth.
In \cite{stable} we proved that there are exactly two tangential
family singularities which are stable under small tangential
deformations (for the Left-Right equivalence relation). 
These singularities, denoted by ${\rm I}$ and ${\rm II}$, 
are represented by $(\xi+t,t^2)$ and $(\xi+t,t^2\xi)$. 
Their envelopes are respectively smooth and have an order $2$
self-tangency. 

Consider the projectivized cotangent bundle $PT^*\R^2$, 
endowed with the standard contact structure and the standard 
Legendre fibration $\pi_L:PT^*\R^2\rightarrow\R^2$. 

\begin{definition*}
The {\it Legendrian graph} of a tangential family is the surface
in $PT^*\R^2$ formed by the Legendrian lifts of the family curves.
\end{definition*}

We remark that the envelope of a tangential family is the 
$\pi_L$-apparent contour of its Legendrian graph.

We say that a Legendrian graph germ is of {\it first type} (resp.,
{\it second type}) if it is generated by a tangential family germ of
the same type, i.e., having a singularity ${\rm I}$
(resp., ${\rm II}$).

\begin{figure}[h]
  \centering
   \scalebox{.3}{\input{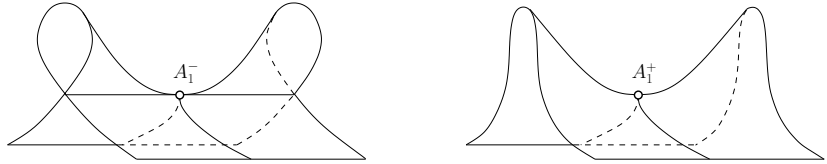}}
  \caption{Double Whitney Umbrellas.}
 \label{fig:1}
\end{figure}

Let us recall that a surface germ has a singularity of type $A^\pm_n$
(resp., $H_n$) if it is diffeomorphic to the surface locally
parameterized by the map germ 
$(\xi,t^2,t^3\pm t\xi^{n+1})$ (resp., by $(\xi,\xi t+t^{3n-1},t^3)$); 
these singularities are simple.
The singularities $A_n^+$ and $A^-_n$ coincide if and only if $n$ is even. 
The singularities $A^\pm_1$, shown in Figure \ref{fig:1}, are called
Double Whitney Umbrellas.

\begin{theorem}\label{thm:1}
The Legendrian graph germs of first type are smooth, while those of second
type have generically a Double Whitney Umbrella singularity $A^\pm_1$.  
The other second type Legendrian graph germs have $A^\pm_n$ or $H_n$ 
singularities for $n\geq 2$ or $n=\infty$.
\end{theorem}

In the statement, ``generically'' means that the second type
Legendrian graph germs for which the claim does not hold form a 
(non-connected) codimension $1$ submanifold in the manifold formed by
the second type Legendrian graph germs.  

\begin{remark*} 
The singularities of map germs
from $\R^2$ to $\R^3$  usually denoted by $B^\pm_n$, $C^\pm_n$, $F_4$
(see \cite{ST2}, \cite{Mond}) appear as singularities of Legendrian graph
generated by non typical tangential families 
(i.e., which are neither of first nor second type). 
For example, the Legendrian graphs of 
$S$-type tangential family germs have $B^\pm_n$ singularities. 
Simple tangential family germs are classified in \cite{simple}. 
\end{remark*}

A Legendrian graph is {\it stable} (under small tangential
deformations) if for every small enough tangential deformation of the
tangential family generating it, the initial and the deformed graphs
are diffeomorphic.  
A similar definition holds for germs. 

\begin{theorem}\label{thm:2}
The Double Whitney Umbrellas $A_1^\pm$ are, in addition to smooth
graphs, the only stable Legendrian graph singularities. 
\end{theorem}
  
We point out that in Mond's general theory of 
maps from the space to the plane \cite{Mond},  
the Double Whitney Umbrellas are not stable. 

A Legendrian graph singularity $L$ is said to be {\it adjacent} to a
Legendrian graph singularity $K$ ($L\rightarrow K$), if
every Legendrian graph in $L$ can
be deformed into a Legendrian graph in $K$ by an arbitrary small
tangential deformation. 
If $L\rightarrow K\rightarrow K'$, the class $L$ is also adjacent to
$K'$. 
In this case we omit the arrow $L\rightarrow K'$. 
The adjacencies of the typical Legendrian graph singularities 
are as follows ($E$ means embedding).
$$
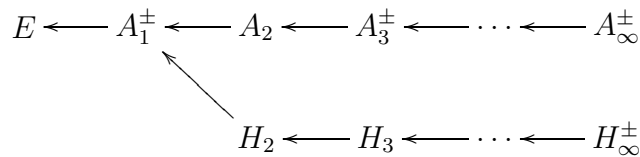
\xymatrix{
 E  &\ar[l] A_1^\pm & \ar[l] A_2    & \ar[l] A_3^\pm  & \ar[l] \dots &
 \ar[l] A_\infty^\pm  \\ 
    &               & \ar[lu] H_2   & \ar[l] H_3      & \ar[l] \dots &
 \ar[l] H_\infty^\pm  }$$

\section{Normal forms of Legendrian graph projections}\label{sct:3}

In this section we study how Legendrian graphs project into the envelopes 
of their generating tangential families. 
In other terms, we find normal forms of typical Legendrian
graphs with respect to an equivalence relation preserving the fiber
structure of $PT^*\R^2$.

\begin{definition*}
The projections of two Legendrian graphs $\Lambda_1$ and $\Lambda_2$ 
by $\pi_L$ are said to be 
{\it equivalent} if there exists a commutative diagram 
$$\begin{CD}
\Lambda_1@>i_1>> PT^*\R^2 @>\pi_L>> \R^2 \\
@VVV     @VVV         @VVV \\
\Lambda_2 @>i_2>> PT^*\R^2 @>\pi_L>> \R^2 \\
\end{CD}$$
in which the vertical arrows are diffeomorphisms and $i_1$, $i_2$ are
inclusions.
\end{definition*}

Such an equivalence is a pair of a diffeomorphism between the
two Legendrian graphs and a diffeomorphism of $PT^*\R^2$ 
fibered over the base $\R^2$ 
(the diffeomorphism is not presumed to be a contactomorphism).  
A similar definition holds for germs. 

Let $\A^*$ be the subgroup of 
$\A:= \Diff (\mathbb{R}^2,0) \times \Diff (\mathbb{R}^3,0)$, 
formed by the pairs $(\phi,\psi)$ such that $\psi$ is fibered with
respect to $\pi$.  
This subgroup inherits the standard action of $\A$ on the maximal
ideal $(\m_{\xi,t})^3$: $(\phi,\psi) \cdot f := \psi\circ f\circ\phi^{-1}$.
Projections of Legendrian graphs are locally equivalent if and
only if their local parameterizations are $\A^*$-equivalent. 

\begin{theorem}\label{thm:3}
The projection germs of the typical Legendrian graphs 
are equivalent to the projection germs of 
the surfaces parameterized by the map germs $f$ in the $3$-space 
$\{x,y,z\}$ by a pencil of lines parallel to the $z$-axis, 
where $f$ is the normal form in the following table, according to the
graph type. 

\begin{center}
\begin{tabular}{c c c c } \hline 
{\rm Type}& {\rm Singularity} & {\rm Normal form}&{\rm Restrictions} \\ \hline 
{\rm I} & {\rm Fold}   & $(\xi,t^2,t)$ & $\emptyset$ \\ 
{\rm II}&$A_1^\pm$       & $(\xi,t^3+t^2\xi+a t\xi^2,t^2+b t^3)$ 
                       &$a\not=-1,0$, $a<1/3$\\ \hline 
\end{tabular}
\end{center}
Moreover, a Legendrian graph germ of second type, parameterized by the above
normal form, has a singularity $A_1^+$ (resp., $A_1^-$) if and
only if $0<a<1/3$ (resp., $-1\not=a<0$).    
\end{theorem}

Typical Legendrian graphs are those having only stable
singularities. 
In Theorem \ref{thm:3},  ``generically'' means that the
second type Legendrian graph germs for which the claim does not hold
form a non-connected codimension $1$ submanifold in
the manifold of all the second type Legendrian graph germs. 

\begin{figure}[h]
  \centering
   \scalebox{.33}{\input{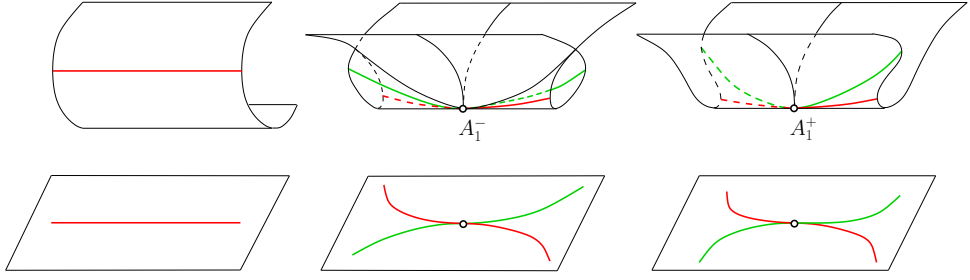}}
        \caption{Typical Legendrian graph projections.}
        \label{fig:2}
\end{figure}

Typical Legendrian graph projections are depicted in
figure \ref{fig:2}.  

\begin{corollary*}
The Fold is the only stable and the only simple 
singularity of Legendrian graphs (with respect to tangential
deformations and Left-Right equivalence relation).
\end{corollary*}

Let us denote by $F_{a,b}$ the $A^\pm_1$ normal form in
Theorem \ref{thm:3} and by $z$ its third coordinate.

\begin{theorem}\label{thm:4}
The map germ $F_{a,b}+ (\mu_1 z,\lambda t+\mu_2 z,0)$ 
is an $\A^*$-miniversal tangential deformation of the normal form $F_{a,b}$, 
provided that $b\not=0$.
\end{theorem}

\begin{remark*}
The above deformation is not the simpler possible, but it has
the property that the parameters $\mu_1$,
$\mu_2$ deform the direction of the projection, but leave fixed the Legendrian
graph, while the parameter $\lambda$ deforms the graph without
changing the projection. 
In particular, the deformation restricted to $\mu_1=\mu_2=0$ provides
an $\A$-miniversal deformation of $F_{a,b}$. 
\end{remark*}

\begin{figure}[h]
  \centering
   \scalebox{.33}{\input{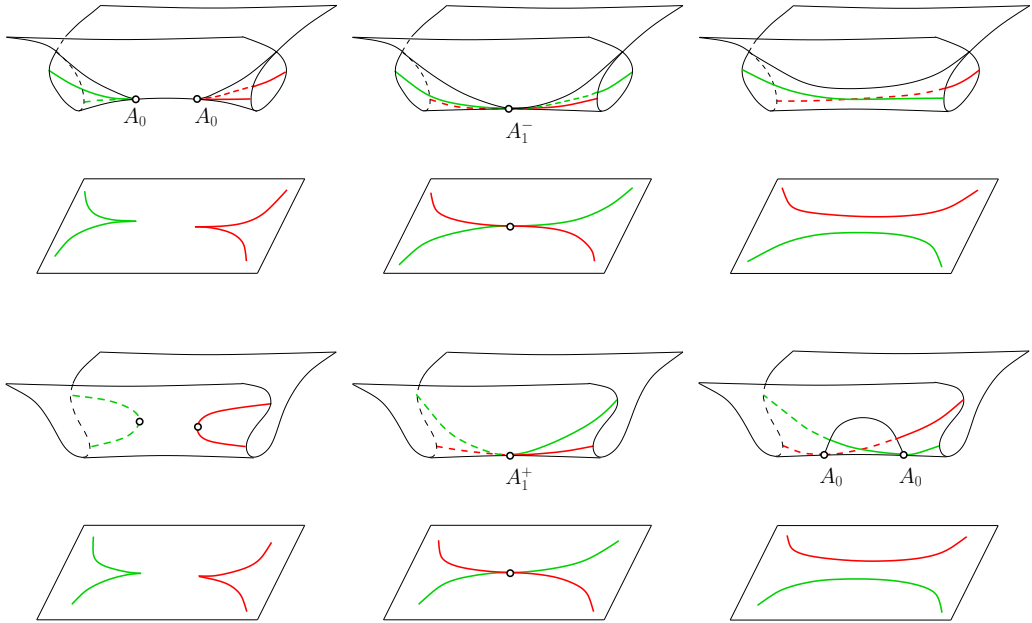}}
  \caption{Legendrian beaks perestroikas.}
 \label{fig:3}
\end{figure}

The second order self-tangency of the envelope of a second
type tangential family germ is not stable under non-tangential 
deformations (see \cite{stable}).  
Under such a deformation, the envelope experiences a beaks  
perestroika, that may be interpreted as the apparent contour in the 
plane of the perestroika occurring to its Legendrian
graph, as shown in figure \ref{fig:3}. 
We call it {\it Legendrian beaks perestroika}. 
Actually, there are two such perestroikas, according to
the sign of $A_1^\pm$.
Figure \ref{fig:3} has been obtained investigating the critical sets 
of the $\A$-miniversal deformation 
$F_{a,b}+(0,\lambda t,0)$ of the projection normal form 
$F_{a,b}$, which leaves unchanged the
direction of the projection. 

\section{Proof of Theorems \ref{thm:1} and \ref{thm:2} }\label{sct:4}
 
We start constructing explicit 
parameterizations of Legendrian graph germs.  

\begin{lemma} 
Every local parameterization of a Legendrian graph is $\A^*$-equivalent
to a map germ of the form  
\begin{equation}\label{eq:1}
\Big ( \xi, k_0 t^2 + (\alpha-k_1)t^3 + k_1 t^2\xi+ \delta(3), 
2k_0 t + (3\alpha-2k_1) t^2 + 
2k_1t\xi+\delta(2) \Big)\ ,
\end{equation}
where $\delta(n)$ denotes any function of $\xi,t$ with zero $n$-jet at
the origin. 
Moreover, the Legendrian graph germ is of first type (resp., of second
type) if and only if if $k_0\not=0$ (resp., $k_0=0$ and
$k_1\not=0,\alpha$).
\end{lemma}

\begin{proof}
Consider a tangential family germ at the origin. 
Up to a coordinate change, preserving the $\A^*$-singularity  of the
graph, we may assume that the family support is locally the $x$-axis. 
For every $x$ small enough, denote by $K(x)$ the curvature at $(x,0)$ 
of the corresponding family curve.  
Now, for $x\rightarrow 0$, let $k_0 +k_1 x+o(x)$ 
and $k_0 x^2+\alpha x^3+o(x^3)$ be the expansions of $K(x)/2$ and of 
the function whose graph (near the origin) 
is the curve associated to $(0,0)$. 

Then, one easily verifies that the Legendrian graph of 
such a tangential family is parameterized by
$(\xi+t,u(\xi,t)+\delta(3), \partial_t u(\xi,t)+\delta(2))$, 
where $u(\xi,t):= k_0 t^2 +\alpha t^3 +   k_1 t^2\xi$. 
This germ can be brought to the required form by
$(\xi,t)\mapsto(\xi-t,t)$. 

Finally, we proved in \cite{stable} that a tangential family
germ is of first type (resp., of second type) 
if and only if $k_0\not=0$ (resp., $k_0=0$ and
$k_1\not=0,\alpha$).   
\end{proof}

We can prove now Theorems \ref{thm:1} and \ref{thm:2}.

\begin{proof}[Proof of Theorem \ref{thm:1}]
If $k_0\not=0$, then the $1$-jet of {\rm (\ref{eq:1})} is $\A$-equivalent to
$(\xi,0,t)$, which is $\A$-sufficient. 
Therefore, Legendrian graph germs of first type are smooth.

We consider now second type Legendrian graph germs (so from now on $k_0=0$). 
First, assume $k_1$ different from the four values 
$0$, $\alpha$, $3\alpha/2$ and $3\alpha$. 
Then {\rm (\ref{eq:1})} is $\A$-equivalent to $(\xi,t^3\pm t\xi^2,t^2)$,
where $\pm$ is the sign of $(k_1-3\alpha)(\alpha-k_1)/k_1^2$. 
Indeed, 
The $3$-jet of (\ref{eq:1}) is $\A$-equivalent to $(\xi, t^3\pm
t\xi^2, t^2)$, which is $\A$-sufficient (see \cite{Mond}, Theorem 1:2).

Hence, the Legendrian graph germs of second type have an 
$A_1^\pm$ singularity whenever  $k_1\not=3\alpha/2,3\alpha$.  
Let us denote by $\hat{\rm II}$ the manifold, formed by all the second
type Legendrian graph local parameterizations. 
The remaining second type graphs belong to the union
of the two submanifolds of $\hat{\rm II}$, defined by $2k_1 = 3\alpha$ and 
$3\alpha = k_1$ (dropping the intersection $\alpha=k_1=0$, whose
elements are not of second type).  
It remains to consider the germs belonging to these two submanifolds. 

If $3\alpha=2k_1\not=0$, the $3$-jet of (\ref{eq:1}) 
is $\A$-equivalent to $(\xi,t^3,t\xi)$.
Then, Mond's classification (see \cite{Mond}, § 4.2.1) implies that 
the map germs (\ref{eq:1}), except those
belonging to an infinite codimension submanifold of $\hat{\rm II}$, 
are $\A$-equivalent to $(\xi, t ^3,t\xi +t^{3n-1})$ for some $n\geq2$. 
On the other hand, when $k_1=3\alpha\not=0$, 
the $2$-jet of (\ref{eq:1}) is $\A$-equivalent to $(\xi,0,t^2)$; 
Mond's classification (see \cite{Mond}, §4.1) implies that the 
map germs of the form {\rm (\ref{eq:1})}, except those
belonging to an infinite codimension submanifold of $\hat{\rm II}$,  
are $\A$-equivalent to $(\xi, t ^3\pm t\xi^{n+1},t^2)$ 
for $n\geq 2$.
\end{proof}

\begin{proof}[Proof of Theorem \ref{thm:2}]
We first show that $A_1^\pm$ singularities are stable. 
It is well known that $(\xi,\lambda t + t^3\pm t^2\xi,t^2)$ is a miniversal 
deformation of $A_1^\pm$ (the singularity being of codimension $1$,
see \cite{Mond}).  
This deformation is not tangential, since it induces a beaks 
perestroika on the corresponding envelope. 
Therefore, every tangential deformation of the singularity is trivial,
due to the envelope stability. 

On the other hand, the non typical Legendrian graphs of second type 
are not stable, due to the adjacencies $A^\pm_{n+1}\rightarrow A^\pm_{n}$, 
$A_{\infty}\rightarrow A^\pm_{n}$, 
$H^\pm_{n+1}\rightarrow H^\pm_{n}$ and 
$H_{\infty}\rightarrow H^\pm_{n}$ (these adjacencies are 
obtained by small tangential deformations). 

Finally, as proven in \cite{stable}, a tangential family germ nether
of first nor second type can be deformed into a second type tangential
family germ via an arbitrary small tangential deformation. 
Hence, its Legendrian graph singularity is adjacent to $A_1^\pm$.  
\end{proof}

\section{Proof of Theorems \ref{thm:3} and \ref{thm:4}}\label{sct:5}

In this section we prove Theorems \ref{thm:3} and \ref{thm:4} 
(for computation details we refer to \cite{these}).  
In order to follow the usual scheme for this reduction, we recall that
a Finite Determinacy Theorem for the $\A^*_3$-equivalence relation 
has been proven by V. V. Goryunov in \cite{goryunov1}; 
this result follows also from Damon's theory 
about {\it nice geometric subgroups of $\A$} (see e.g. \cite{Damon1}) . 

A map germ $f\in(\m_{\xi,t})^3$ defines, by $f^*g:=g\circ f$, a
homomorphism from the ring $\E_{x,y,z}$ of the function germs in the
target to the ring $\E_{\xi,t}$ of the function germs in the source. 
Hence, every $\E_{\xi,t}$-module has a structure of
$\E_{x,y,z}$-module via this homomorphism. 
We define the {\it extended tangent space} of $f$ as usual by    
$$T_e\A^*(f):= \langle \partial_\xi f,\partial_t
f\rangle_{\E_{\xi,t}}  + 
f^*(\E_{x,y})\times f^*(\E_{x,y}) \times f^*(\E_{x,y,z}) \ . $$ 
Note that $T_e\A^*(f)$ is an $\E_{x,y}$-module, being in general
neither an $\E_{\xi,t}$-module nor an $\E_{x,y,z}$-module. 
The {\it reduced tangent space} $T_r\A^*(f)$ of $f$ is 
by definition the $\mathscr E_{x,y}$-submodule
of $T_e\A^*(f)$ defined by   
$\mathfrak g_+  + \mathcal M^*$, 
where $\mathfrak g_+$ is the space of all the vector field
germs having positive order (see \cite{AGV} for definitions) 
and $\mathcal M^*$ is the following $\E_{x,y}$-module: 
$$f^*\left(\m^2_{x,y}\oplus \langle y \rangle_\R\right) \times 
f^*\left(\m^2_{x,y}\oplus \langle x \rangle_\R\right) \times 
f^*\left(\m^2_{x,y,z}\oplus\langle x,y \rangle_\R\right) \ .$$

The main tool in the proof of Theorems
\ref{thm:3} and \ref{thm:4} is the following easy fact.

\begin{lemma}\label{lemma:2}
Consider $f\in(\m_{\xi,t})^3$ and $R$ a triple of homogeneous 
polynomials of degree $p$, $q$ and $r$, such that $R\in
T_r\A^*(f)/(\m_{\xi,t}^{p+1}\times \m_{\xi,t}^{q+1}\times\m_{\xi,t}^{r+1})$.   
Then the $(p,q,r)$-jets of $f$ and $f+R$ are  $\A^*$-equivalent. 
\end{lemma}

We can start now the proof of Theorem \ref{thm:3}. 

\begin{proof}[Proof of Theorem \ref{thm:3}.]
We consider first Legendrian graphs of first type tangential family
germs ($k_0\not=0$).  
Then the $2$-jet of (\ref{eq:1}) is $\A^*$-equivalent to 
$(\xi,t^2,t)$, which is $\A^*$-sufficient, since its reduced tangent
space contains $\m_{\xi,t}^2\times \m^3_{\xi,t}\times \m^2_{\xi,t}$
(Lemma \ref{lemma:2}). 

We consider now Legendrian graphs of second type tangential family
germs ($k_0=0$).   
In this case, every map germ (\ref{eq:1}) is $\A^*$-equivalent to 
$\left(\xi,t^3+t^2\xi+at\xi^2 + \delta(3), t^2+\delta(2)\right)$, 
where $a:=(\alpha-k_1)(k_1-3\alpha)/k_1^2$. 
We remark that $a<1/3$; indeed, we have   
$1-3a =(3\alpha-2k_1)^2/k_1^2>0$, since $3\alpha\not=2k_1$. 
Actually, a further computation shows that its $(2,4,3)$-jet is 
$\A^*$-equivalent to $F_{a,b}$, for a suitable $b\in\R$ 
($F_{a,b}$ is the normal for defined in Section \ref{sct:3}).  
Hence, the statement follows from Lemma \ref{lemma:2} and the next
inclusion, which holds for  $a\not=-1,0,1/3$: 
\begin{equation}\label{eq:3}
\m_{\xi,t}^3\times\m_{\xi,t}^5\times\m_{\xi,t}^4\subset T_r\A^*(F_{a,b}) \ . 
\end{equation}
When the Legendrian graph has an  $A_1^\pm$ singularity, the conditions
$a\not=0,1/3$ are automatically fulfilled. 
On the other hand, $a\not=-1$ is a new condition, equivalent to
$\alpha=0$, giving rise to the submanifold for which Theorem \ref{thm:3} does not
hold.  
\end{proof}

\begin{proof}[Proof of Theorem \ref{thm:4}.] 
Since $\langle \xi, t^3+t^2\xi+a t \xi^2\rangle_{\E_{\xi,t}} =\langle
\xi,t^3\rangle_{\E_{\xi,t}} $ and 
$$\langle \xi, t^3+t^2\xi+a t \xi^2,t^2+bt^3 \rangle_{\E_{\xi,t}} =\langle
\xi,t^2\rangle_{\E_{\xi,t}}\ ,$$ the well known Preparation Theorem of
Mather--Malgrange (see e.g. \cite{AGV}) implies that 
$\E_{\xi,t}$ is generated by $\{t,t^2\}$ as
$\E_{x,y}$-module and by $t$ as  $\E_{x,y,z}$-module. 
Hence, we have: 
$$\E_{\xi,t}^3= F^*_{a,b}(\E_{x,y})\cdot \left\{
\left(\begin{smallmatrix} t\\ 0\\ 0 \end{smallmatrix}\right) , 
\left(\begin{smallmatrix} 0\\ t\\ 0 \end{smallmatrix}\right) , 
\left(\begin{smallmatrix} t^2\\ 0\\ 0 \end{smallmatrix}\right) ,
\left(\begin{smallmatrix} 0\\ t^2\\ 0 \end{smallmatrix}\right) 
\right\} +
F^*_{a,b}(\E_{x,y,z})\cdot \left\{
\left(\begin{smallmatrix} 0\\ 0\\ t \end{smallmatrix}\right) 
\right\} \ . $$  
For $b\not=0$, we obtain 
$$\E_{\xi,t}^3= T_e\A^*(F_{a,b}) \oplus \R\cdot \left\{
\left(\begin{smallmatrix} 0\\ t\\ 0 \end{smallmatrix}\right) , 
\left(\begin{smallmatrix} z\\ 0\\ 0 \end{smallmatrix}\right) , 
\left(\begin{smallmatrix} 0\\ z\\ 0 \end{smallmatrix}\right) 
\right\}  \ . $$
This proves the Theorem.
\end{proof}



\end{document}